\documentclass{article}
\usepackage{amsfonts}
\usepackage{amsmath}

\newtheorem{theorem}{Theorem}

\newtheorem{proposition}[theorem]{Proposition}
\newtheorem{remark}[theorem]{Remark}

\input{tcilatex}
\begin{document}

\title{A Special Nonlinear Connection in Second Order Geometry}
\author{Nicoleta BR\^{I}NZEI \\
%EndAName
"Transilvania" University, Brasov, Romania\\
e-mail: nico.brinzei@rdslink.ro}
\maketitle


\section{Introduction}

As shown in \cite{Rahula}, nonlinear connections on bundles can be a
powerful tool for integrating systems of differential equations. A way of
obtaining them is that of deriving them out of the respective systems of
DE's, in particular, out of variational principles, \cite{Anto}, \cite{Miron}%
, \cite{MironSapporo}. In particular an ODE system of order 2 on a manifold $%
M$ induces a nonlinear connection on its tangent bundle. A remarkable
example is here the Cartan nonlinear connection of a Finsler space, which
has the property that its autoparallel curves correspond to geodesics of the
base manifold:%
\begin{equation*}
\dfrac{\delta y^{i}}{dt}:=\dfrac{dy^{i}}{dt}+N_{~j}^{i}y^{j}=0.
\end{equation*}

Further, an ODE system of order three determines a nonlinear connection on
the second order tangent (jet) bundle $T^{2}M=J_{0}^{2}(\mathbb{R},M)$. For
instance, Craig-Synge equations (R. Miron, \cite{Miron})%
\begin{equation*}
\dfrac{d^{3}x^{i}}{dt^{3}}+3!G^{i}(x,\dot{x},\ddot{x})=0,
\end{equation*}%
lead to:

\qquad a) Miron's connection:%
\begin{equation}
\underset{\left( 1\right) }{M}\underset{}{_{j}^{i}}=\dfrac{\partial G^{i}}{%
\partial y^{\left( 2\right) j}},\,\underset{\left( 2\right) }{M}\underset{}{%
_{j}^{i}}=\frac{1}{2}\left( S\underset{\left( 1\right) }{M}\underset{}{%
_{j}^{i}}+\underset{\left( 1\right) }{M}\underset{}{_{m}^{i}}\underset{%
\left( 1\right) }{M}\underset{}{_{j}^{m}}\right) ,  \label{Mir}
\end{equation}%
where $S=y^{\left( 1\right) i}\dfrac{\partial }{\partial x^{i}}+2y^{\left(
2\right) i}\dfrac{\partial }{\partial y^{\left( 1\right) i}}-3G^{i}\dfrac{%
\partial }{\partial y^{\left( 2\right) i}}.$

\qquad b) Buc\u{a}taru's connection%
\begin{equation*}
\underset{1}{M}\overset{}{_{\;j}^{i}}=\frac{\partial G^{i}}{\partial y^{(2)j}%
},\underset{2}{M}\overset{}{_{\;j}^{i}}=\frac{\partial G^{i}}{\partial
y^{(1)j}}.
\end{equation*}%
W.r.t. the last one, if$\ G^{i}$ are the coefficients of a spray on $T^{2}M$
(i.e., 3-homogeneous functions), then the Craig-Synge equations can be
interpreted as:%
\begin{equation}
\dfrac{\delta y^{(2)i}}{dt}=0,  \label{*}
\end{equation}%
where $\dfrac{\delta y^{(2)i}}{dt}:=\dfrac{dy^{(2)i}}{dt}+~\underset{1}{M}%
\overset{}{_{\;j}^{i}}\dfrac{dy^{(1)j}}{dt}+~\underset{2}{M}\overset{}{%
_{\;j}^{i}}\dfrac{dx^{j}}{dt}.$

In Miron's and Buc\u{a}taru's approaches, nonlinear connections on $T^{2}M,$
are obtained from a Lagrangian of order 2, $L(x,\dot{x},\ddot{x}),$ by
computing the first variation of its integral of action.

Here, we propose a different approach, which, we consider, is at least as
interesting as the above one from the point of view of Mechanics - just by
having in view that the big majority of known Lagrangians are of order one.

Namely, we start with a first order Lagrangian $L(x,\dot{x})$ and compute
its second variation; out of it, we obtain a nonlinear connection on $%
T^{2}M, $ such that the obtained distributions correspond to extremal curves
and their deivations, respectively.

As a remark, our nonlinear connection is also suitable for modelling the
solutions of a (globally defined) ODE system, not necessarily attached to a
certain Lagrangian, together with their deviations.

\bigskip

More precisely, in the following our aims are:

\begin{enumerate}
\item to obtain the Jacobi equations for the trajectories%
\begin{equation*}
\dfrac{\delta y^{i}}{dt}=F^{i}(x,y)
\end{equation*}%
(for extremal curves of a 2-homogeneous Lagrangian $L(x,\dot{x})$ in
presence of external forces).

\item (main result): to build a nonlinear connection such that:

\begin{equation*}
w\in \mathcal{X}(M)\text{ Jacobi field along }c\Leftrightarrow ~\dfrac{%
\delta w^{(2)i}}{dt}=0,
\end{equation*}

where $\dfrac{d}{dt}$ denotes directional derivative w.r.t. $\dot{c}$ and $%
\dfrac{\delta w^{(2)i}}{dt}=\dfrac{1}{2}\dfrac{d^{2}w^{i}}{dt^{2}}+~\underset%
{1}{M}\overset{}{_{\;j}^{i}}\dfrac{dw^{i}}{dt}+~\underset{2}{M}\overset{}{%
_{\;j}^{i}}w^{j}.$
\end{enumerate}

\textbf{Properties:}

\textbf{I. }For Finsler spaces $(M,F)$, $c$ is a geodesic of $M$ if and only
if its extension $T^{2}M$ is horizontal.

\textbf{II. }For a vector field $w$ along a geodesic $c$ on $M,$ we have:

\begin{enumerate}
\item $\dfrac{\delta w^{i}}{dt}=0,$ if and only if $w$ is parallel along $%
\dot{c}=y.$

\item $\dfrac{\delta w^{(2)i}}{dt}=0\ $if and only if $w$ is a Jacobi field
along $c.$
\end{enumerate}

\section{Tangent Bundle and 2-Tangent Bundle}

\qquad Let $M$ be a real differentiable manifold of dimension $n$ and class $%
\mathcal{C}^{\infty }$; the coordinates of a point $x\in M$ in a local chart 
$\left( U,\phi \right) $ will be denoted by $\phi \left( x\right) =\left(
x^{i}\right) ,$ $i=1,...,n.$ Let $TM$ be its tangent bundle and $%
(x^{i},y^{i})$ the coordinates of a point in a local chart.

The \emph{2-tangent bundle }(or \emph{2-osculator bundle})\emph{\ }$%
(T^{2}M,\pi ^{2},M)$ is the space of jets of order two in a fixed point $%
t_{0}\in I\subset \mathbb{R},$ of functions $f:I\rightarrow M,$ $t\mapsto
(f^{i}(t)),$ (\cite{12}-\cite{17}).

In a local chart, a point $p$ of $T^{2}M$ will have the coordinates $%
(x^{i},y^{i},y^{(2)i}).$ This is,%
\begin{equation*}
x^{i}=f^{i}(t_{0}),~\ \ y^{i}=\dot{f}^{i}(t_{0}),~\ \ y^{(2)i}=\dfrac{1}{2}%
\overset{\cdot \cdot }{f}\overset{}{^{i}}(t_{0}).
\end{equation*}%
Obviously, $\left( T^{2}M,\pi ^{2},M\right) $ is a differentiable manifold
of class $\mathcal{C}^{\infty }$ and dimension $3n,$ and $TM=Osc^{1}M$ can
be identified with a submanifold of $T^{2}M$.

For a curve $c:[0,1]\rightarrow M,$ $t\mapsto (x^{i}(t))$, let us denote:

\begin{itemize}
\item by $\widehat{c}$ its \textit{extension} to the tangent bundle $TM:$%
\begin{equation*}
\widehat{c}:[0,1]\rightarrow M,t\mapsto (x^{i}(t),\dot{x}^{i}(t));
\end{equation*}%
along $\widehat{c}$, there holds:%
\begin{equation*}
y^{i}=\dot{x}^{i}(t);
\end{equation*}

\item by $\widetilde{c}$ its \textit{extension to }$T^{2}M$:%
\begin{equation*}
\widetilde{c}:[0,1]\rightarrow T^{2}M,~\ \ \ t\mapsto (x^{i}(t),\dot{x}%
^{i}(t),\dfrac{1}{2}\overset{\cdot \cdot }{x}\overset{}{^{i}}(t));
\end{equation*}%
along such an extension curve, there holds%
\begin{equation*}
y^{i}(t)=\dot{x}^{i}(t),~\ \ y^{(2)i}(t)=\dfrac{1}{2}\overset{\cdot \cdot }{x%
}\overset{}{^{i}}(t).
\end{equation*}%
\bigskip
\end{itemize}

\section{Nonlinear connections on $TM$}

A \textit{nonlinear (Ehresmann) connection on }$TM$ is a splitting of the
tangent space $T(TM)$ in every point $p\in TM$ into a direct sum%
\begin{equation}
T_{p}(TM)=N\left( p\right) \oplus V\left( p\right) ,
\end{equation}%
each one of dimension $n.$ This generates two distributions:

\begin{itemize}
\item the \textit{horizontal distribution }$N:p\mapsto N(p);$

\item the \textit{vertical} \textit{distribution }$V:p\mapsto V(p).$
\end{itemize}

A local adapted basis to the above decomposition is: 
\begin{equation*}
B=\{\frac{\delta }{\delta x^{i}},\frac{\partial }{\partial y^{i}}\},
\end{equation*}%
where:%
\begin{equation}
\begin{array}{l}
\dfrac{\delta }{\delta x^{i}}=\dfrac{\partial }{\partial x^{i}}-N_{~i}^{j}%
\dfrac{\partial }{\partial y^{j}}%
\end{array}%
\end{equation}%
With respect to changes of local coordinates on $M,$ $\dfrac{\delta }{\delta
x^{i}}$ transform by the same rule as vector fields on $M:$ $\dfrac{\delta }{%
\delta x^{i}}=\dfrac{\partial \widetilde{x}^{j}}{\partial x^{i}}\dfrac{%
\delta }{\delta \widetilde{x}^{j}}.$

The dual basis of $B$ is $B^{\ast }=\left\{ dx^{i},\delta y^{i}\right\} $,
given by%
\begin{equation}
\delta y^{i}=dy^{i}+N_{~j}^{i}dx^{j}.
\end{equation}

The quantities $N_{~j}^{i}$ are called the \textit{coefficients }of the
nonlinear connection $N.$

Any vector field $X\in \mathcal{X}\left( TM\right) $ is represented in the
local adapted basis as 
\begin{equation}
X=X^{\left( 0\right) i}\dfrac{\delta }{\delta x^{i}}+X^{\left( 1\right) i}%
\frac{\partial }{\partial y^{i}}.
\end{equation}

Similarly, a 1-form $\omega \in \mathcal{X}^{\ast }\left( TM\right) $ will
be decomposed as%
\begin{equation}
\omega =\omega _{i}^{\left( 0\right) }dx^{i}+\omega _{i}^{\left( 1\right)
}\delta y^{\left( 1\right) i}.
\end{equation}

In particular, if $\widehat{c}:t\rightarrow (x^{i}(t),y^{i}(t))$ is an
extension curve to $TM$, then its tangent vector field is%
\begin{equation}
T=\dfrac{dx^{i}}{dt}\delta _{(0)i}+\dfrac{\delta y^{(1)i}}{dt}\delta _{(1)i}.
\end{equation}

We should mention an important result, (R. Miron, \cite{MironSapporo}):

\begin{proposition}
Let $L=L(x,\dot{x})$ be a nondegenerate Lagrangian: $\det (\dfrac{\partial
^{2}L}{\partial y^{i}\partial y^{j}})\not=0,$ and $g_{ij}=\dfrac{1}{2}\dfrac{%
\partial ^{2}L}{\partial y^{i}\partial y^{j}},$ the induced (Lagrange)
metric tensor. Then, the equations of evolution of a mechanical system with
the Lagrangian $L$ and the external force field $F=F_{i}(x,\dot{x})dx^{i}$
are%
\begin{equation}
\dfrac{d^{2}x^{i}}{dt^{2}}+2G^{i}(x,\dot{x})=\dfrac{1}{2}F^{i}(x,\dot{x}),
\label{evolution}
\end{equation}%
where 
\begin{equation*}
2G^{i}=\dfrac{1}{2}g^{is}(\dfrac{\partial ^{2}L}{\partial y^{s}\partial x^{j}%
}y^{j}-\dfrac{\partial L}{\partial x^{j}}),
\end{equation*}%
is the canonical semispray of the Lagrange space $(M,L)$ and $%
F^{i}=g^{ij}F_{j}.$
\end{proposition}

\textbf{(!) }In the following, we shall use the above result in the case
when $G$ is a spray, i.e., its components $G^{i}$ are 2-homogeneous w.r.t. $%
y^{i}:$%
\begin{equation*}
2G^{i}=\dfrac{\partial G^{i}}{\partial y^{j}}y^{j}.
\end{equation*}%
Then, \cite{Anto}, \cite{Shen}, \cite{Mi-An}, the quantities $N_{~j}^{i}=%
\dfrac{\partial G^{i}}{\partial y^{j}}$ are the coeficients of a nonlinear
connection on $TM,$ and, moreover, the equations (\ref{evolution}) take the
form:%
\begin{equation}
\dfrac{\delta y^{i}}{dt}=\dfrac{1}{2}F^{i}.  \label{evolution1}
\end{equation}

In particular, if there are no external forces, $F^{i}=0,$ then the extremal
curves $t\mapsto x^{i}(t)$ of the Lagrangian $L$ have horizontal extensions
and vice-versa: horizontal extension curves $\widehat{c}$ project onto
solution curves of the Euler-Lagrange equations of $L.$

\section{Nonlinear connections on $T^{2}M$}

A \textit{nonlinear connection on }$T^{2}M$ is a splitting of the tangent
space in every point $p$ into three subspaces:%
\begin{equation}
T_{p}(T^{2}M)=N_{0}\left( p\right) \oplus N_{1}\left( p\right) \oplus
V_{2}\left( p\right) ,  \label{Np}
\end{equation}%
each one of dimension $n.$ This generates three distributions:

\begin{itemize}
\item the \textit{horizontal distribution }$N_{0}:p\mapsto N(p);$

\item the $v_{1}$-\textit{distribution }$N_{1}:p\mapsto N_{1}(p);$

\item the $v_{2}$-\textit{distribution }$V_{2}:p\mapsto V_{2}(p).$
\end{itemize}

We denote by $h=v_{0},$ $v_{1}$ and $v_{2}$ the projectors corresponding to
the above distributions

Let $\mathcal{B}$ denote a local adapted basis to the decomposition (\ref{Np}%
): 
\begin{equation*}
\mathcal{B}=\{\delta _{(0)i}:=\frac{\delta }{\delta x^{i}}=\frac{\delta }{%
\delta y^{\left( 0\right) i}},\delta _{(1)i}:=\frac{\delta }{\delta y^{i}}%
,\delta _{(2)i}:=\frac{\delta }{\delta y^{\left( 2\right) i}}\},
\end{equation*}%
this is, $N_{0}=Span(\delta _{(0)i}),$ $N_{1}=Span(\delta _{(1)i}),$ $%
V_{2}=Span(\delta _{(2)i}).$ The elements of the adapted basis are locally
expressed as%
\begin{equation}
\left\{ 
\begin{array}{l}
\delta _{(0)i}=\dfrac{\delta }{\delta x^{i}}=\dfrac{\partial }{\partial x^{i}%
}-\underset{\left( 1\right) }{N}\underset{}{_{i}^{j}}\dfrac{\partial }{%
\partial y^{j}}-\underset{\left( 2\right) }{N}\underset{}{_{i}^{j}}\dfrac{%
\partial }{\partial y^{\left( 2\right) j}} \\ 
\delta _{(1)i}=\dfrac{\delta }{\delta y^{\left( 1\right) i}}=~\ \ \ \ \ \ \
\ \ \ \ \ \ \ \dfrac{\partial }{\partial y^{i}}-\underset{\left( 1\right) }{N%
}\underset{}{_{i}^{j}}\dfrac{\partial }{\partial y^{\left( 2\right) j}} \\ 
\delta _{(2)i}=\dfrac{\delta }{\delta y^{\left( 2\right) i}}=~\ \ \ \ \ \ \
\ \ \ \ \ \ \ \ \ \ \ \ \ \ \ \ \ \ \ \ \ \ \ \ \ \ \ \dfrac{\partial }{%
\partial y^{\left( 2\right) i}}.%
\end{array}%
\right.
\end{equation}%
With respect to changes of local coordinates on $M,$ $\delta _{(\alpha )i},$ 
$\alpha =0,1,2,$ transform by the same rule as vector fields on $M:$ $\delta
_{(\alpha )i}=\dfrac{\partial \widetilde{x}^{j}}{\partial x^{i}}\widetilde{%
\delta }_{(\alpha )j}.$

The dual basis of $\mathcal{B}$ is $\mathcal{B}^{\ast }=\left\{
dx^{i},\delta y^{\left( 1\right) i},\delta y^{\left( 2\right) i}\right\} $,
given by%
\begin{eqnarray}
\delta y^{\left( 0\right) i} &=&dx^{i}, \\
\delta y^{\left( 1\right) i} &:&=\delta y^{i}=dy^{i}+\underset{\left(
1\right) }{M}\underset{}{_{j}^{i}}dx^{j}, \\
\delta y^{\left( 2\right) i} &=&dy^{\left( 2\right) i}+\underset{\left(
1\right) }{M}\underset{}{_{j}^{i}}dy^{\left( 1\right) j}+\underset{\left(
2\right) }{M}\underset{}{_{j}^{i}}dx^{j}.
\end{eqnarray}

Also, w.r.t. local chart changes, $\delta y^{\left( \alpha \right) i}$
behave exactly as 1-forms on $M:$ $\delta y^{\left( \alpha \right) i}=\dfrac{%
\partial x^{i}}{\partial \widetilde{x}^{j}}\widetilde{\delta }y^{(\alpha
)j}. $

The quantities $\underset{\left( 1\right) }{N}\underset{}{_{i}^{j}},$ $%
\underset{\left( 2\right) }{N}\underset{}{_{i}^{j}}$ are called the \textit{%
coefficients }of the nonlinear connection $N,$ while $\underset{\left(
1\right) }{M}\underset{}{_{j}^{i}}$ and $\underset{\left( 2\right) }{M}%
\underset{}{_{j}^{i}}$ are called its \textit{dual} coefficients.

Then, a vector field $X\in \mathcal{X}\left( T^{2}M\right) $ is represented
in the local adapted basis as 
\begin{equation}
X=X^{\left( 0\right) i}\delta _{(0)i}+X^{\left( 1\right) i}\delta
_{(1)i}+X^{\left( 2\right) i}\delta _{(2)i},
\end{equation}%
with the three right terms (called \textit{d-vector fields}) belonging to
the distributions $N,$ $N_{1}$ and $V_{2}$ respectively.

A 1-form $\omega \in \mathcal{X}^{\ast }\left( T^{2}M\right) $ will be
decomposed as%
\begin{equation}
\omega =\omega _{i}^{\left( 0\right) }dx^{i}+\omega _{i}^{\left( 1\right)
}\delta y^{i}+\omega _{i}^{\left( 2\right) }\delta y^{\left( 2\right) i}.
\end{equation}%
Similarly, a tensor field $T\in \mathcal{T}_{s}^{r}\left( T^{2}M\right) $
can be split \ with respect to (\ref{Np}) into components , named \textit{%
d-tensor fields}, which transform w.r.t. local coordinate changes in the
same way as tensors on $M.$

In particular, if $\widetilde{c}:t\rightarrow
(x^{i}(t),y^{i}(t),y^{(2)i}(t)) $ is an extension curve, then its tangent
vector field is%
\begin{equation}
T=\dfrac{dx^{i}}{dt}\delta _{(0)i}+\dfrac{\delta y^{(1)i}}{dt}\delta _{(1)i}+%
\dfrac{\delta y^{(2)i}}{dt}\delta _{(2)i}.  \label{T}
\end{equation}

Our aim is to give a precise meaning to the equality $v_{2}(\widetilde{c}%
)=0. $

\section{Berwald linear connection}

Let $G^{i}$ be the coefficients of a spray on $TM,$ and%
\begin{equation*}
N_{~j}^{i}=\dfrac{\partial G^{i}}{\partial y^{j}},
\end{equation*}%
the coefficients of the induced nonlinear connection (on $TM$).

Let also%
\begin{equation*}
L_{~jk}^{i}(x,y)=\dfrac{\partial N_{~j}^{i}}{\partial y^{k}}=\dfrac{\partial
G^{i}}{\partial y^{j}\partial y^{k}},
\end{equation*}%
the local coefficients of the induced Berwald linear connection on $TM.$

Now, let on $T^{2}M,$ $\underset{\left( 1\right) }{N}\underset{}{_{\,\,j}^{i}%
}~=~N_{~j}^{i}(x,y^{(1)}),$ as above. The \textit{Berwald connection }on $%
T^{2}M$ , \cite{Bucataru1}, is the linear connection defined by%
\begin{equation}
\left\{ 
\begin{array}{c}
D_{\delta _{(0)k}}\delta _{(\alpha )j}=L_{~jk}^{i}\delta _{(\alpha )i}, \\ 
D_{\delta _{(\beta )k}}\delta _{(\alpha )j}=0.%
\end{array}%
\right.
\end{equation}%
This is, with the notations in \cite{Miron}, the coefficients of the Berwald
linear connection are $B\Gamma (N)=(L_{~jk}^{i},0,0).$

For extension curves $\widetilde{c},$ we can express the $v_{1}$ component $%
\dfrac{\delta y^{i}}{dt}$ (the \textit{geometrical acceleration, \cite{Lewis}%
}) by means of the Berwald covariant derivative:%
\begin{equation*}
\dfrac{Dy^{i}}{dt}:=D_{\overset{\cdot }{\widetilde{c}}}y^{i}=\dfrac{\delta
y^{i}}{dt}.
\end{equation*}

Let $\mathbb{T}$ denote its torsion tensor, and:%
\begin{equation*}
R_{~jk}^{i}=v_{1}\mathbb{T}(\delta _{(0)k},\delta _{(0)j})=\delta
_{(0)k}N_{~j}^{i}-\delta _{(0)j}N_{~k}^{i};
\end{equation*}%
also, let $\mathbb{R}$ be the curvature tensor, and%
\begin{eqnarray*}
R_{j~kl}^{~i}\delta _{(0)i} &=&h\mathbb{R}(\delta _{(0)l},\delta
_{(0)k})=L_{~jk}^{i}-\delta
_{k}L_{~jl}^{i}+L_{~jk}^{m}L_{~ml}^{i}-L_{~jl}^{m}L_{~mk}^{i}, \\
P_{j~kl}^{~i}\delta _{(0)i}~ &=&h\mathbb{R}(\delta _{(1)l},\delta
_{(0)k})=\delta _{(1)l}L_{~jk}^{i}.
\end{eqnarray*}%
its local components, which will be used in the following.

\section{Jacobi equations for systems with external forces}

Let us consider the following PDE system%
\begin{equation}
A^{i}\equiv \dfrac{\partial ^{2}X^{i}}{\partial x^{k}\partial x^{j}}%
y^{j}y^{k}+\{2\delta _{j}^{i}y^{(2)k}+a_{~j}^{i}(x,y,y^{(2)})y^{k}\}\dfrac{%
\partial X^{j}}{\partial x^{k}}+b_{~j}^{i}(x,y,y^{(2)})X^{j}=0,
\label{syst3}
\end{equation}%
where:

\begin{itemize}
\item the unknown functions $X^{i}=X^{i}(x)$ are the local coordinates of a
vector field on $M$;

\item w.r.t. local coordinate changes, $A^{i}$ trasnsform as: $\widetilde{A}%
^{~j}=\dfrac{\partial \widetilde{x}^{j}}{\partial x^{i}}A^{j}$ (i.e., $A^{i}$
are the components of a d-vector field on $T^{2}M$).
\end{itemize}

Then, it is easy to check

\begin{proposition}
The quantities 
\begin{eqnarray*}
\underset{\left( 1\right) }{M}\underset{}{_{\,\,j}^{i}}(x,y,y^{(2)}) &=&%
\dfrac{1}{2}a_{~j}^{i}(x,y,y^{(2)}), \\
\underset{\left( 2\right) }{M}\underset{}{_{\,\,j}^{i}}(x,y,y^{(2)}) &=&%
\dfrac{1}{2}b_{~j}^{i}(x,y,y^{(2)})
\end{eqnarray*}%
are the dual coefficients of a nonlinear connection on $T^{2}M.$
\end{proposition}

Now, let us interpret (\ref{syst3}). If $c:t\in \lbrack 0,1]\rightarrow M$
is a curve, $\widetilde{c}$ is its extension to $T^{2}M,$ and $X$ is a
vector field along $c,$ then (\ref{syst3}) can be reexpressed in the
following way:%
\begin{equation*}
\dfrac{d^{2}X^{i}}{dt^{2}}+~2\underset{\left( 1\right) }{M}\underset{}{%
_{\,\,j}^{i}}\dfrac{dX^{j}}{dt}+~2\underset{\left( 2\right) }{M}\underset{}{%
_{\,\,j}^{i}}X^{j}=0.
\end{equation*}%
It is now convenient to denote%
\begin{equation*}
X^{(2)j}=\dfrac{1}{2}\dfrac{dX^{j}}{dt};
\end{equation*}%
with these notation, (\ref{syst3}) can be brought to the much more familiar
form%
\begin{equation}
\dfrac{\delta X^{(2)i}}{dt}=0,  \label{syst3'}
\end{equation}%
which leads to the following interpretation: the integral curves of $X$ have
extensions of second order with vanishing $v_{2}$-components.

\bigskip

Let us suppose that we know \textit{a priori} a nonlinear connection on $TM,$
with 1-homogeneous coefficients $N_{~j}^{i}=\dfrac{\partial G^{i}}{\partial
y^{j}}.$ Let $c:[0,1]\rightarrow M,$ $t\mapsto x^{i}(t)$ be a curve and $%
\widehat{c}:[0,1]\rightarrow TM,$ $t\mapsto (x^{i}(t),y^{i}(t)=\dot{x}%
^{i}(t))$ its extension to $TM.$ Let us suppose that $x^{i}$ are solutions
for the system of ODE's%
\begin{equation}
\dfrac{\delta y^{i}}{dt}=F^{i}(x,y),  \label{Lagrange}
\end{equation}%
where $F^{i}$ are the components of a d-vector field on $M$ (for commodity
of notations, we omitted the $\dfrac{1}{2}$ in front of $F$)

Let $x^{i}(t,u)$ be a variation of $c$ (not necessarily with fixed
endpoints) $y^{i}=y^{i}=\dfrac{\partial x^{i}}{\partial t},$ and%
\begin{equation*}
w^{i}(t)=\dfrac{\partial x^{i}}{\partial u}|_{u=0}
\end{equation*}%
its associated deviation vector field.

Then, by commuting partial derivatives, we have as immediate consequences%
\begin{equation*}
\dfrac{\partial y^{i}}{\partial u}=\dfrac{dw^{i}}{dt}=:2w^{(2)i},~\ \dfrac{%
\partial ^{2}y^{i}}{\partial t^{2}}=\dfrac{d^{2}w^{i}}{dt^{2}}=2\dfrac{%
dw^{(2)i}}{dt}~\ \ etc.
\end{equation*}

Let us denote%
\begin{eqnarray*}
\dfrac{\delta y^{i}}{\partial t} &=&\dfrac{\partial y^{i}}{\partial t}+~M%
\underset{}{_{\,\,j}^{i}}(x,y)y^{j}, \\
\dfrac{\delta y^{i}}{\partial u} &=&\dfrac{\partial y^{i}}{\partial u}+~M%
\underset{}{_{\,\,j}^{i}}(x,y)w^{j}, \\
\dfrac{\delta w^{i}}{\partial t} &=&\dfrac{\partial w^{i}}{\partial t}+~M%
\underset{}{_{\,\,j}^{i}}(x,y)w^{j};
\end{eqnarray*}%
then, $\dfrac{\delta y^{i}}{\partial t},$ $\dfrac{\delta y^{i}}{\partial u}$
and $\dfrac{\delta w^{i}}{\partial t}$ define d-vector fields on $T^{2}M$
(the covariant derivatives "with reference vector $\dot{c}~$" of $\dot{c}$
and $w$),

It is also immediate that the last two covariant derivatives coincide:%
\begin{equation*}
\dfrac{\delta y^{i}}{\partial u}=\dfrac{\delta w^{i}}{\partial t},
\end{equation*}%
which can be interpreted in terms of Berwald connection on $TM$ as%
\begin{equation*}
D_{\tfrac{\partial \widetilde{c}}{\partial u}}(h(\frac{\partial \widetilde{c}%
}{\partial t}))=D_{\tfrac{\partial \widetilde{c}}{\partial t}}(h(\frac{%
\partial \widetilde{c}}{\partial u})).
\end{equation*}%
By applying $D_{\tfrac{\partial \widetilde{c}}{\partial t}}$ again to the
above quality and by permuting covariant derivatives, we get

\begin{theorem}
The components of the deviation vector fields $w^{i}$ of the trajectories%
\begin{equation}
\dfrac{\delta y^{i}}{dt}=F^{i}(x,y),
\end{equation}%
satisfy the Jacobi-type equation%
\begin{equation}
\dfrac{D^{2}w^{i}}{dt^{2}}=\dfrac{DF^{i}}{\partial u}%
|_{u=0}+y^{h}y^{j}R_{h~jk}^{~i}w^{k}.  \label{Jacobi}
\end{equation}
\end{theorem}

It now remains to express (\ref{Jacobi}) in natural coordinates and
"collect" the coefficients of $w^{i}$ and its derivatives, in order to
obtain a nonlinear connection on $T^{2}M.$

By a direct computation, we obtain

\begin{proposition}
The deviation vector fields of the trajectories $\dfrac{\delta y^{i}}{dt}%
=F^{i}(x,y)$ are solutions of ODE system:%
\begin{equation*}
\begin{array}{l}
\dfrac{d^{2}w^{i}}{dt^{2}}+(2M_{~j}^{i}-\dfrac{\partial F^{i}}{\partial y^{j}%
})\dfrac{dw^{j}}{dt}+ \\ 
+(\mathbb{C}%
(M_{~j}^{i})+M_{~k}^{i}M_{~j}^{k}-y^{h}y^{k}R_{h~jk}^{~i}+L_{~kj}^{i}F^{j}-%
\dfrac{\partial F^{i}}{\partial x^{j}})w^{j}=0,%
\end{array}%
\end{equation*}%
where $\mathbb{C}=y^{k}\dfrac{\partial }{\partial x^{k}}+2y^{(2)k}\dfrac{%
\partial }{\partial y^{k}}$
\end{proposition}

By re-expressing the above in terms of the nonlinear connection, $%
R_{~hjk}^{i}y^{h}=R_{~jk}^{i},$ $L_{~kj}^{i}=\dfrac{\partial M_{~k}^{i}}{%
\partial y^{j}},$ we get

\begin{theorem}
\begin{enumerate}
\item The quantities%
\begin{eqnarray*}
\underset{\left( 1\right) }{M}\underset{}{_{\,\,j}^{i}}~ &=&\dfrac{1}{2}%
(2M_{~j}^{i}-\dfrac{\partial F^{i}}{\partial y^{j}}), \\
\underset{\left( 2\right) }{M}\underset{}{_{\,\,j}^{i}}~ &=&\dfrac{1}{2}(%
\mathbb{C}(M_{~j}^{i})+M_{~k}^{i}M_{~j}^{k}-y^{k}R_{~jk}^{~i}+\dfrac{%
\partial M_{~k}^{i}}{\partial y^{j}}F^{j}-\dfrac{\partial F^{i}}{\partial
x^{j}})
\end{eqnarray*}%
are the dual coefficients of a nonlinear connection on $T^{2}M.$
\end{enumerate}
\end{theorem}

\begin{theorem}
With respect to this nonlinear connection, the extensions of Jacobi fields
attached to (\ref{Lagrange}) have vanishing $v_{2}$ components:%
\begin{equation*}
\dfrac{1}{2}\dfrac{d^{2}w^{i}}{dt^{2}}+~\underset{\left( 1\right) }{M}%
\underset{}{_{\,\,j}^{i}}\dfrac{dw^{j}}{dt}+~\underset{\left( 2\right) }{M}%
\underset{}{_{\,\,j}^{i}}w^{j}=0.
\end{equation*}
\end{theorem}

By denoting%
\begin{equation*}
w^{(1)i}=\dfrac{dw^{i}}{dt},~\ \ \ w^{(2)i}=\dfrac{1}{2}\dfrac{d^{2}w^{i}}{%
dt^{2}},
\end{equation*}%
the above relations can be reexpressed as:

\begin{equation*}
\dfrac{\delta w^{(2)i}}{dt}=0.
\end{equation*}

\section{Deviations of geodesics}

Let $N_{~j}^{i}=\dfrac{\partial G^{i}}{\partial y^{j}}$ be a nonlinear
connection on $TM,$ coming from a spray.

If $F=0,$ then we deal with deviations of autoparallel curves,%
\begin{equation*}
\dfrac{\delta y^{i}}{dt}=0.
\end{equation*}%
In this case,%
\begin{eqnarray*}
\underset{\left( 1\right) }{M}\underset{}{_{\,\,j}^{i}}~ &=&M_{~j}^{i}, \\
\underset{\left( 2\right) }{M}\underset{}{_{\,\,j}^{i}}~ &=&\dfrac{1}{2}(%
\mathbb{C}(M_{~j}^{i})+M_{~k}^{i}M_{~j}^{k}-y^{j}R_{~jk}^{i}),
\end{eqnarray*}%
and we notice that our nonlinear connection differs only by the term $%
-y^{j}R_{~jk}^{i}$ from Miron's one, \cite{Miron}.

\begin{remark}
Along an extension curve $\widetilde{c},$ there hold the equalities%
\begin{equation*}
\dfrac{\delta y^{i}}{dt}=\dfrac{Dy^{i}}{dt},~\ \ \dfrac{\delta y^{(2)i}}{dt}=%
\dfrac{D^{2}y^{i}}{dt^{2}},
\end{equation*}%
where $\dfrac{D}{dt}$ denotes the covariant derivative associated to the
Berwald connection.
\end{remark}

\begin{remark}
Also, for a vector field $w$ along the projection $c$ on $M,$ we have%
\begin{equation*}
\dfrac{\delta w^{i}}{dt}=\dfrac{Dw^{i}}{dt}.
\end{equation*}
\end{remark}

\textbf{Conclusions:}

\begin{enumerate}
\item $c$ is a geodesic if and only if its extension $T^{2}M$ is horizontal.
\end{enumerate}

For a vector field $w$ along a geodesic $c$ on $M,$ we have:

\begin{enumerate}
\item $\dfrac{\delta w^{i}}{dt}=0,$ if and only if $w$ is parallel along $%
\dot{c}=y.$

\item $\dfrac{\delta w^{(2)i}}{dt}=0\ $if and only if $w$ is a Jacobi field
along $c.$
\end{enumerate}

This is, the extension to $T^{2}M$ of a Jacobi field $w$ on $M$ will have
the form%
\begin{equation*}
W=w^{i}\delta _{(0)i}+\dfrac{\delta w^{i}}{dt}\delta _{(1)i}.
\end{equation*}

\section{External forces in locally Minkovsian spaces}

Let $(M,L(y))$ be a locally Minkovski space.

Then, $M_{~j}^{i}=0,$ $L_{~jk}^{i}=0$ (for the Berwald connection), \cite%
{Anto}, \cite{Shen}. In presence of an external force field, the evolution
equations of a mechanical system will take the form%
\begin{equation}
\dfrac{dy^{i}}{dt}=F^{i}(x,y).  \label{Mink}
\end{equation}

In this case, with the above notations, we have%
\begin{eqnarray*}
\underset{\left( 1\right) }{M}\underset{}{_{\,\,j}^{i}}~ &=&-\dfrac{1}{2}%
\dfrac{\partial F^{i}}{\partial y^{j}}, \\
\underset{\left( 2\right) }{M}\underset{}{_{\,\,j}^{i}}~ &=&-\dfrac{1}{2}%
\dfrac{\partial F^{i}}{\partial x^{j}},
\end{eqnarray*}%
and the deviations of trajectories (\ref{Mink}) are, simply:%
\begin{equation*}
\dfrac{d^{2}w^{i}}{dt^{2}}-~\dfrac{\partial F^{i}}{\partial y^{j}}\dfrac{%
dw^{j}}{dt}-~\dfrac{\partial F^{i}}{\partial x^{j}}w^{j}=0.
\end{equation*}

\end{document}